\pdfoutput=1
\RequirePackage{ifpdf}
\ifpdf 
\documentclass[pdftex]{sigma}
\else
\documentclass{sigma}
\fi

\numberwithin{equation}{section}

\newtheorem{Theorem}{Theorem}[section]
\newtheorem*{Theorem*}{Theorem}

 { \theoremstyle{definition}
\newtheorem{Definition}[Theorem]{Definition}

 }

\newcommand\mylabel[1]{\label{#1}}
\newcommand\thm[1]{\ref{thm:#1}}

\newcommand\eqn[1]{(\ref{eq:#1})}

\newcommand\refdef[1]{\ref{def:#1}}

\begin{document}


\renewcommand{\thefootnote}{}

\newcommand{\arXivNumber}{2406.12049}

\renewcommand{\PaperNumber}{097}

\FirstPageHeading

\ShortArticleName{Combinatorial Interpretations of Cranks of Overpartitions}

\ArticleName{Combinatorial Interpretations of Cranks\\ of Overpartitions and Partitions\\ without Repeated Odd Parts\footnote{This paper is a~contribution to the Special Issue on Basic Hypergeometric Series Associated with Root Systems and Applications in honor of Stephen C.~Milne's 75th birthday. The~full collection is available at \href{https://www.emis.de/journals/SIGMA/Milne.html}{https://www.emis.de/journals/SIGMA/Milne.html}}}

\Author{Frank G. GARVAN and Rishabh SARMA}

\AuthorNameForHeading{F.G.~Garvan and R.~Sarma}

\Address{Department of Mathematics, University of Florida,\\
P.O. Box 118105, Gainesville, FL 32611-8105, USA}
\Email{\href{mailto:fgarvan@ufl.edu}{fgarvan@ufl.edu}, \href{mailto:rishabh.sarma@ufl.edu}{rishabh.sarma@ufl.edu}}
\URLaddress{\url{https://qseries.org/fgarvan/}, \url{https://people.clas.ufl.edu/rishabh-sarma/}}

\ArticleDates{Received June 19, 2024, in final form October 12, 2024; Published online October 26, 2024}

\Abstract{We give combinatorial interpretations of two residual cranks of overpartitions defined by Bringmann, Lovejoy and Osburn in 2009 analogous to the crank of partitions given by Andrews and the first author in 1988. As a consequence, we give new versions of their definitions without adjusted weights. Furthermore, we investigate the combinatorial interpretation of an $M_2$-crank of partitions without repeated odd parts and explore connections of these statistics with their companion rank counterparts and the tenth order mock theta functions of Ramanujan.}


\Classification{05A15; 05A17; 05A30; 11P81}

\begin{flushright}
\begin{minipage}{81mm}
\it
Dedicated to our friend and colleague Steve Milne\\ on the occasion of his 75th
birthday
\end{minipage}
\end{flushright}

\section{Introduction}\label{section1}

A partition of a positive integer $n$ is a finite non-increasing sequence of positive integers $\lambda_1 \geq \lambda_2 \geq \cdots \geq \lambda_k$ such that $\lambda_1 + \lambda_2 + \cdots + \lambda_k = n$. We denote $p(n)$ to be the number of (unrestricted) partitions of $n$. For example, $p(5) = 7$, where the seven partitions of $5$ are
\begin{align*}
5,\quad 4+1,\quad 3+2,\quad 3+1+1,\quad 2+2+1,\quad 2+1+1+1,\quad 1+1+1+1+1.
\end{align*}
In 2004, Corteel and Lovejoy \cite{Co-Lo} introduced the concept of overpartitions. An overpartition is a partition in which the first occurrence of a number may be overlined.
For example, the $14$ overpartitions of $4$ are
\begin{gather*}
4,\quad  \overline{4},\quad  3 + 1,\quad  \overline{3} + 1,\quad  3 + \overline{1},\quad  \overline{3} + \overline{1},\quad  2 + 2,\quad  \overline{2} + 2,\quad  2 + 1 + 1,
\\
\overline{2} + 1 + 1,\quad  2 + \overline{1} + 1,\quad  \overline{2} + \overline{1} + 1,\quad 1 + 1 + 1 + 1,\quad  \overline{1} + 1 + 1 + 1.
\end{gather*}

In 2009, Bringmann, Lovejoy and Osburn \cite{Br-Lo-Os} defined what they called ``residual cranks'' for overpartitions and studied the quasimodularity properties of the moments of these two cranks. In their paper, the authors mention that the crank for partitions arose because of its relation to Ramanujan's congruences and the non-existence of such simple congruences for overpartitions was studied by Choi \cite{Ch09}.

The elusive crank for partitions was found by Andrews and the first author in search of a~combinatorial explanation of the mod $11$ congruence of Ramanujan \cite{An-Ga88}. The crank of a partition is defined as the largest part if there are no ones in the partition and otherwise the number of parts larger than the number of ones minus the number of ones. It divides the partitions of $5n+4$, $6n+5$ and $11n+6$ into $5$, $7$ and $11$ equinumerous classes respectively, thus providing an almost combinatorial explanation of all three partition congruences of Ramanujan. Let $M(m,n)$ denote the number of partitions of $n$ with crank $m$. Many authors mistakenly write
\begin{equation*}
\sum_{n=0}^\infty \sum_{m=-n}^n
M(m,n) z^m q^n =
\dfrac{(q;q)_{\infty}}{(zq;q)_{\infty}\bigl(z^{-1}q;q\bigr)_{\infty}},
\end{equation*}
where we use the standard $q$-notation
\begin{align*}
(a;q)_{\infty}= \prod_{k=0}^\infty \bigl(1-aq^k\bigr),
\qquad
(a;q)_{n}= \frac{(a;q)_{\infty}}{(aq^n;q)_{\infty}}.
\end{align*}
This is a mistake since
\begin{align*}
&{}\sum_{n=0}^\infty
\sum_{m=-n}^n N_V(m,n)z^mq^n=
\dfrac{(q;q)_{\infty}}{(zq;q)_{\infty}\bigl(z^{-1}q;q\bigr)_{\infty}}
=1+\bigl(-1+z+z^{-1}\bigr)q+\cdots,
\end{align*}
where $N_V(m, n)$ is the number of vector partitions of $n$ with vector
crank $m$ (see \cite{An-Ga88,Ga86,Ga88}).
Note that
\begin{equation*}
N_V(-1,1)=1,\qquad N_V(0,1)=-1,\qquad N_V(1,1)=1,
\end{equation*}
but
\begin{equation*}
M(-1,1)=1,\qquad M(0,1)=0,\qquad M(1,1)=0.
\end{equation*}
In \cite{Ga88b}, it was shown that $N_V(m,n)\geq 0$ for $n>1$. The proof of this inequality led to the solution of Dyson's crank conjecture in terms of partitions and the definition of the crank \cite{An-Ga88} so that
\begin{equation*}
N_V(m,n)=M(m,n)
\end{equation*}
for $n>1$.

We define
\begin{equation*}
C(z,q) :=\dfrac{(q;q)_{\infty}}{(zq;q)_{\infty}\bigl(z^{-1}q;q\bigr)_{\infty}}
=1+\bigl(-1+z+z^{-1}\bigr)q+\sum_{n=2}^\infty \sum_{m=-n}^n M(m,n)z^mq^n.
\end{equation*}
We note that the function $C(z,q)$ appears in Ramanujan's lost notebook (see \cite{Ga86,Ga88b,Ga88}). The~lost notebook contains the $5$-dissection of $C(\zeta_5,q)$ \big(where $\zeta_5=\exp\bigl(\frac{2\pi {\rm i}}{5}\bigr)$\big) and the analogue for Dyson's rank function. The rank of a partition is the largest part minus the number of parts. In the lost notebook, these functions are also connected with Ramanujan's fifth order mock theta functions, via rank and crank differences. For example, if $N(a,m,n)$ and $M(a,m,n)$ denote the number of partitions of $n$ with rank and respectively crank congruent to $a$ modulo~$m$, then
\begin{align*}
&{}\sum_{n\geq 0}3(N(1,5,5n)-N(2,5,5n))q^n-\sum_{n\geq 0}(M(0,5,5n)-M(1,5,5n))q^n
=\chi_0(q)-2,
\end{align*}
where
\begin{equation*}
\chi_0(q)=\sum\limits_{n\geq 0}\frac{q^n}{\bigl(q^{n+1};q\bigr)_{n}},
\end{equation*}
is one of Ramanujan's fifth order mock theta functions. This result is equivalent to Ramanujan's equation~(4.11) on \cite[p.~36]{Ga88b}. This is the first mock theta conjecture \cite{An-Ga89} which was proved by Hickerson \cite{Hi}.

Bringmann, Lovejoy and Osburn \cite{Br-Lo-Os} were in fact the first to define such a notion of crank for overpartitions from which they then deduce congruence properties for combinatorial functions which can be expressed in terms of the second overpartition rank moment and the corresponding residual crank moment. Several authors have since considered these residual cranks and have worked on their generalizations, finding and proving inequalities between the moments of these functions in conjunction with other overpartition statistics, among other problems.

The coefficients $\overline{M}(m,n)$ and $\overline{M2}(m,n)$ are defined by the following:
\begin{align*}
&{}\overline{C}(z,q) = \sum_{n=0}^\infty \sum_m \overline{M}(m,n) z^m q^n=(-q;q)_{\infty}C(z;q)=\dfrac{\bigl(q^2;q^2\bigr)_{\infty}}{(zq;q)_{\infty}\bigl(z^{-1}q;q\bigr)_{\infty}},
\\
&{}\overline{C2}(z,q) = \sum_{n=0}^\infty \sum_m \overline{M2}(m,n) z^m q^n=\frac{(-q;q)_{\infty}}{\bigl(q;q^2\bigr)_{\infty}}C\bigl(z;q^2\bigr)=\dfrac{(-q;q)_{\infty}\bigl(q^2;q^2\bigr)_{\infty}}{\bigl(q;q^2\bigr)_{\infty}\bigl(zq^2;q^2\bigr)_{\infty}\bigl(z^{-1}q^2;q^2\bigr)_{\infty}}.
\end{align*}
Bringmann, Lovejoy and Osburn \cite{Br-Lo-Os} describe $\overline{M}(m,n)$ (resp.\ $\overline{M2}(m,n)$) as the number of overpartitions of $n$ with first (resp.\ second) residual crank equal to $m$. They claim that the first residual crank of an overpartition is obtained by taking the crank of the subpartition consisting of the non-overlined parts. Similarly, the second residual crank is obtained by taking the crank of the subpartition consisting of all of the even nonoverlined parts divided by two. The authors state that they make the appropriate modifications based on the fact that for partitions we have $M(0,1)=-1$ and $M(-1,1)=M(1,1)=1$. However, what they claim is not counting overpartitions but involves a modified count since
\begin{equation*}
M(m,n)\neq N_V(m,n)
\end{equation*}
for $(m,n)=(0,1),(1,1)$. In this paper, we give combinatorial interpretations of $\overline{M}(m,n)$ and $\overline{M2}(m,n)$ solely in terms of overpartitions with no adjusted weights. These interpretations are given in Theorems \thm{M1} and \thm{M2} below. Our new definitions of the first and second residual cranks depend on the overlined parts of the overpartition if there are no non-overlined parts.

We describe Bringmann, Lovejoy and Osburn's modified count at the end of Section~\ref{section2}.

For an ordinary partition $\pi$, let $\ell(\pi)$ denote the largest part of $\pi$, $\omega(\pi)$
denote the number of ones in $\pi$, and $\mu(\pi)$ denote the number of parts of $\pi$
larger than $\omega(\pi)$.
Then
\begin{equation*}
\operatorname{crank}(\pi)=\begin{cases}
\ell(\pi) & \mbox{if $\omega(\pi)=0$}, \\
\mu(\pi)-\omega(\pi) & \mbox{otherwise}.
\end{cases}
\end{equation*}
Let $\mathcal{P}$, $\mathcal{PD}$, $\mathcal{PO}$ and $\mathcal{PE}$ denote the sets of unrestricted partitions, partitions into distinct parts, partitions into odd parts and partitions into even parts, respectively. Let $\varnothing$ denote the empty partition.
\begin{Definition}
Let $\pi$ be a partition where all parts are distinct. Define $\lambda(\pi)$ as
\begin{equation*}
\lambda(\pi)=\begin{cases}
0 & \mbox{if $\pi=\varnothing$ or $\ell(\pi), \,\ell(\pi)-1$ are both parts},
\\
1 & \mbox{if $\pi=1$ or has only one part or otherwise.}
\end{cases}
\end{equation*}
Then writing an overpartition as an element of $\vec{\pi}=(\pi_1,\pi_2) \in \mathcal{PD} \times \mathcal{P}$, the unmodified first residual crank of overpartitions is given by
\begin{align}
\mylabel{eq:crank1}
\operatorname{crank1}(\vec{\pi})=\begin{cases}
\operatorname{crank}(\pi_2) & \mbox{if $\pi_2 \neq \varnothing$},
\\
\lambda(\pi_1) & \mbox{if $\pi_2 = \varnothing$.}
\end{cases}
\end{align}
\end{Definition}
An example calculating the first residual crank $\operatorname{crank1}(\vec{\pi})$ for overpartitions $\vec{\pi}$ of $3$ is given towards the end of Section~\ref{section2}.
\begin{Theorem}
\mylabel{thm:M1}
The number of overpartitions $\vec{\pi}$ of $n$ with $\operatorname{crank1}(\vec{\pi})=m$ is $\overline{M}(m,n)$ for all $n$.
\end{Theorem}

\begin{Definition}
\mylabel{def:crank2}
Let $\pi$ be a partition, where all parts are distinct. We define $\kappa(\pi)$ into the following subcases.

Case I: $\ell(\pi) \geq 4$,
\begin{equation*}
\kappa(\pi)=\begin{cases}
0 & \mbox{if $\ell(\pi)-1$ or $\ell(\pi)-2$ are parts of $\pi$}, \\
1 & \mbox{otherwise.}
\end{cases}
\end{equation*}

Case II: $\ell(\pi) = 3$ or $2$,
\begin{equation*}
\kappa(\pi)=\begin{cases}
1 & \mbox{if there is only one part},
\\
0 & \mbox{otherwise.}
\end{cases}
\end{equation*}

Case III: $\ell(\pi) = 1$, then $\kappa(\pi)=0$.

Then, writing an overpartition as an element of $\vec{\pi}=(\pi_1,\pi_2,\pi_3) \in \mathcal{PD} \times \mathcal{PE} \times \mathcal{PO}$, the unmodified second residual crank of overpartitions is given by
\begin{align}
\mylabel{eq:crank2}
\operatorname{crank2}(\vec{\pi})=\begin{cases}
\operatorname{crank}\bigl(\frac{\pi_2}{2}\bigr) & \mbox{if $\pi_2 \neq \varnothing$},
\\
\kappa(\pi_1) & \mbox{if $\pi_2 = \varnothing$  and $\pi_1 \neq \varnothing$,}
\\
0 & \mbox{otherwise},
\end{cases}
\end{align}
where a partition $\frac{\pi}{2}$ is obtained by dividing each part of the partition $\pi \in \mathcal{PE}$ by $2$.
\end{Definition}
An example calculating the second residual crank $\operatorname{crank2}(\vec{\pi})$ for overpartitions $\vec{\pi}$ of $4$ is given at the end of Section~\ref{section3}.

\begin{Theorem}
\mylabel{thm:M2}
The number of overpartitions $\vec{\pi}$ of $n$ with $\operatorname{crank2}(\vec{\pi})=m$ is $\overline{M2}(m,n)$ for all~$n$.
\end{Theorem}

A notable significance of the first residual crank of overpartitions lies in the fact that it arises in the formulas for overpartition rank differences. Some instances of overpartition rank differences are mock theta functions of order $10$ as noted by Lovejoy and Osburn \cite[p.~197]{Lo-Os10}:
\begin{align*}
&\sum_{n\geq 0}\bigl(\overline{N}(0,5,5n+1)-\overline{N}(2,5,5n+1)\bigr)q^n=2\phi(q),
\\
&\sum_{n\geq 0}\bigl(\overline{N}(0,5,5n+4)+\overline{N}(1,5,5n+4)-2\overline{N}(2,5,5n+4)\bigr)q^{n+1}=2\psi(q),
\end{align*}
where
\begin{equation*}
\phi(q)=\sum\limits_{n\geq 0}\frac{q^{\binom{n+1}{2}}}{\bigl(q;q^2\bigr)_{n+1}}
\qquad\mbox{and}\qquad
\psi(q)=\sum\limits_{n\geq 0}\frac{q^{\binom{n+2}{2}}}{\bigl(q;q^2\bigr)_{n+1}}
\end{equation*}
are two tenth order mock theta functions of Ramanujan.
Similarly, some instances of overpartition rank differences can be written as a sum of a first residual crank difference and a mock theta function of order $10$ analogous to the fact that some partition rank differences can be written as a sum of the Andrews--Garvan crank differences and a mock theta function of order $5$. For example, if $\overline{N}(a,m,n)$ and $\overline{M}(a,m,n)$ denote the number of overpartitions of $n$ with rank and respectively first residual crank congruent to $a$ modulo $m$, then one can deduce that
\begin{gather*}
\sum_{n\geq 0}\bigl(\overline{N}(0,5,5n+1)-\overline{N}(1,5,5n+1)\bigr)q^n-\sum_{n\geq 0}\bigl(\overline{M}(0,5,5n+1)-\overline{M}(1,5,5n+1)\bigr)q^n \\
\qquad{}=3\phi(q).
\end{gather*}
The interested reader may try to see that this follows from the overpartition rank difference identities of Lovejoy and Osburn {\cite[equations~(1.6) and~(1.11)]{Lo-Os08}}, overpartition first residual crank difference identities due to the first author and Jennings-Shaffer \cite[Theorem~2.10]{Ga-Sh} and the identity
\begin{equation*}
\phi(q)=a_2(q)+2qh\bigl(q^2,q^5\bigr)
\end{equation*}
for the tenth order mock theta function $\phi(q)$ due to Choi (see \cite[p. 534]{Ch99} for details).

Similarly, one can deduce that
\begin{align*}
&{}\sum_{n\geq 0}\bigl(\overline{N}(1,5,5n+4)-\overline{N}(2,5,5n+4)\bigr)q^n-\sum_{n\geq 0}\bigl(\overline{M}(0,5,5n+4)-\overline{M}(2,5,5n+4)\bigr)q^n \\
&{}\qquad{}=3q^{-1}\psi(q).
\end{align*}
This follows from the overpartition rank difference identities of Lovejoy and Osburn \cite[equations~(1.6) and~(1.11)]{Lo-Os08}, overpartition first residual crank difference identities due to Garvan and Jennings-Shaffer \cite[Theorem~2.10]{Ga-Sh} and the identity
\begin{equation*}
\psi(q)=a_1(q)+2qh\bigl(q,q^5\bigr)
\end{equation*}
for the tenth order mock theta function $\psi(q)$ due to Choi (see \cite[p.~533]{Ch99} for details).

The two other tenth order mock theta functions of Ramanujan have connections with the $M_2$-rank and crank of partitions without repeated odd parts which we explore in a later section.

On the other hand, the second residual crank is the companion statistic to the $M_2$-rank of overpartitions. More specifically, several instances of linear combinations of $M_2$-rank differences for overpartitions can be written in terms of linear combinations of second residual crank differences for overpartitions and many such combinations have nice Jacobi product representations. A systematic study of these identities will be taken up in a later paper. Further, the second residual crank of overpartitions also arises in the formula for the two variable overpartition rank generating function. Specifically, replacing $z$ by a root of unity, the two variable overpartition rank generating function can be written as the sum of a function due to Zwegers, a Mordell integral of a theta function and the second residual crank generating function with appropriate multipliers (see \cite[Corollary 2.2]{JS} for details).

The rest of our paper is organized into three sections as follows. In Sections~\ref{section2} and~\ref{section3}, following the lines of argument of Andrews and the first author \cite{An-Ga88}, we deduce combinatorial interpretations of the first and second residual crank of overpartitions respectively by dissecting their generating functions to arrive at the definitions and theorems stated above in the introduction. We close each of these two sections with an illustrative example. In Section~\ref{section4}, we define an $M_2$-crank of partitions without repeated odd parts and explore its combinatorial interpretation and relations with other two tenth order mock theta functions of Ramanujan and the $M_2$-rank of partitions without repeated odd parts studied by Lovejoy and Osburn \cite{Lo-Os09}.

\section{Proof of Theorem \thm{M1}}\label{section2}

\begin{proof}[\unskip\nopunct]
Cauchy's identity \cite[p.\ 17]{An98} states
\begin{equation*}
\sum\limits_{n\geq 0}\frac{(a;q)_n}{(q;q)_{n}}t^n=\frac{(at;q)_{\infty}}{(t;q)_{\infty}}.
\end{equation*}
Using this, we have
\begin{align*}
C(z,q) &{}=\dfrac{(q;q)_{\infty}}{(zq;q)_{\infty}\bigl(z^{-1}q;q\bigr)_{\infty}}
\\
&{}=\frac{1-q}{(zq;q)_{\infty}}\frac{\bigl(q^2;q\bigr)_{\infty}}{\bigl(z^{-1}q;q\bigr)_{\infty}}
\\
&{}=\frac{1-q}{(zq;q)_{\infty}}\sum\limits_{n=0}^{\infty}\frac{(zq;q)_n}{(q;q)_n}\bigl(z^{-1}q\bigr)^n
\\
&{}=\frac{1-q}{(zq;q)_{\infty}}+\sum\limits_{n=1}^{\infty}\frac{z^{-n}q^n}{\bigl(q^2;q\bigr)_{n-1}\bigl(zq^{n+1};q\bigr)_{\infty}}
\\
&{}=(1-q)\sum\limits_{n=0}^{\infty}\frac{z^{n}q^n}{(q;q)_{n}}+\sum\limits_{n=1}^{\infty}\frac{z^{-n}q^n}{\bigl(q^2;q\bigr)_{n-1}\bigl(zq^{n+1};q\bigr)_{\infty}}
\\
&{}=(1-q)+zq+\sum\limits_{n=2}^{\infty}\frac{z^{n}q^n}{\bigl(q^2;q\bigr)_{n-1}}+\sum\limits_{n=1}^{\infty}\frac{z^{-n}q^n}{\bigl(q^2;q\bigr)_{n-1}\bigl(zq^{n+1};q\bigr)_{\infty}}.
\end{align*}
Thus,
\begin{align}
\mylabel{eq:crankexp}
C(z,q)=(1-q+zq)+\sum_{n=1}^\infty \sum_{m=-n}^n M(m,n) z^m q^n.
\end{align}
\\
Since $\overline{C}(z,q)=(-q;q)_{\infty}C(z;q)$, we look for a combinatorial interpretation of $(-q;q)_{\infty}(1-q+zq)$. We know
\begin{equation*}
(-q;q)_{\infty}=1+\sum\limits_{m=1}^{\infty}(1+q)\cdots\bigl(1+q^{m-1}\bigr)q^m,
\end{equation*}
where $m$ corresponds to the largest part of a partition. Therefore,
\begin{align}
\mylabel{eq:distexp}
&\sum\limits_{m=1}^{\infty}(1+q)\cdots\bigl(1+q^{m-1}\bigr)q^m=(-q;q)_{\infty}-1.
\end{align}
Removing the first part, we get
\begin{align}
\mylabel{eq:RFP1}
&\sum\limits_{m=2}^{\infty}(1+q)\cdots\bigl(1+q^{m-2}\bigr)\bigl(1+q^{m-1}\bigr)q^m=(-q;q)_{\infty}-1-q.
\end{align}
Also, we have the following identity:
\begin{align}
\mylabel{eq:EID1}
&\sum\limits_{m=2}^{\infty}(1+q)\cdots\bigl(1+q^{m-2}\bigr)q^m=q\sum\limits_{m=2}^{\infty}(1+q)\cdots\bigl(1+q^{m-2}\bigr)q^{m-1}=q((-q;q)_{\infty}-1).
\end{align}
Then, using equations \eqn{distexp}--\eqn{EID1}, we get
\begin{align*}
(-q;q)_{\infty}(1-q+zq)&{}=1+zq+zq ((-q;q)_{\infty}-1)+((-q;q)_{\infty}-1-q)-q((-q;q)_{\infty}-1)
\\
&{}=1+zq+\sum\limits_{m=2}^{\infty}(1+q)\cdots\bigl(1+q^{m-2}\bigr)q^m\bigl(z+q^{m-1}\bigr)
\\
&{}=\sum\limits_{\pi \in \mathcal{PD}}z^{\lambda(\pi)}q^{|\pi|},
\end{align*}
where $|\pi|$ is the sum of parts of $\pi$. We write an overpartition as an element of $\vec{\pi}=(\pi_1,\pi_2) \in \mathcal{PD} \times \mathcal{P}$. Thus
\begin{align*}
\overline{C}(z,q)&{}=\sum_{n=0}^\infty \sum_{m=-n}^n M(m,n) z^m q^n
\\
&{}=(-q;q)_{\infty}(1-q+zq)+(-q;q)_{\infty}\sum_{n=1}^\infty
\sum_{m=-n}^n M(m,n) z^m q^n
\\
&{}=\sum\limits_{\pi_1 \in \mathcal{PD}}z^{\lambda(\pi_1)}q^{|\pi_1|}+\sum\limits_{\pi_1 \in \mathcal{PD}}q^{|\pi_1|}\sum\limits_{\substack{\pi_2 \in \mathcal{P}\\ \pi_2 \neq \varnothing}}z^{\operatorname{crank}(\pi_2)}q^{|\pi_2|}
\\
&{}=\sum\limits_{\vec{\pi} = (\pi_1,\pi_2) \in \mathcal{PD}\times \mathcal{P}}z^{\operatorname{crank1}(\vec{\pi})}q^{|\vec{\pi}|},
\end{align*}
where we have employed equation \eqn{crank1} and $|\vec{\pi}|=|(\pi_1,\pi_2)|=|\pi_1|+|\pi_2|.$ Theorem \thm{M1} follows.
\end{proof}

We present an example.
\begin{table}[hbt!]\centering
\caption{First residual crank of overpartitions of $3$.}\label{t1}
\vspace{2mm}
\begin{tabular}{ c c c c c c }
Overpartitions of $3\,(\vec{\pi})$ & $\pi_1$ & $\lambda(\pi_1)$ & $\pi_2$ & $\operatorname{crank}(\pi_2)$ & $\operatorname{crank1}(\vec{\pi})$ \\
 $3$                         & $\varnothing$ &     &  $3$   & $3$ & $3$\\
 $\overline{3}$              & $3$    & $1$ & $\varnothing$ &     & $1$\\
 $2+1$                       & $\varnothing$ &     & $2+1$  & $0$ & $0$\\
 $\overline{2}+1$            & $2$    &     &  $1$   & \hspace{-3mm}$-1$& \hspace{-3mm}$-1$\\
 $2+\overline{1}$            & $1$    &     &  $2$   & $2$ & $2$\\
 $\overline{2}+\overline{1}$ & $2+1$  & $0$ & $\varnothing$ &     & $0$\\
 $1+1+1$                     & $\varnothing$ &     & $1+1+1$& \hspace{-3mm}$-3$& \hspace{-3mm}$-3$\\
 $\overline{1}+1+1$          & $1$    &     & $1+1$  & \hspace{-3mm}$-2$& \hspace{-3mm}$-2$
\end{tabular}
\end{table}

From Table \ref{t1}, we find
\begin{alignat*}{4}
&\overline{M}(-3,3)=1,\qquad&&\overline{M}(-2,3)=1,\qquad&&\overline{M}(-1,3)=1,&
\\
&\overline{M}(0,3)=2,&
\\
&\overline{M}(1,3)=1,\qquad&&\overline{M}(2,3)=1,\qquad&&\overline{M}(3,3)=1.&
\end{alignat*}
This agrees with Bringmann, Lovejoy and Osburn's modified count which we now describe. Since
\begin{equation*}
N_V(-1,1)=1,\qquad N_V(0,1)=-1,\qquad N_V(1,1)=M(1,1)=1,
\end{equation*}
if an overpartition of $n$ has one $1$ as an unoverlined part, then that overpartition contributes a~$-1$ to the count of $\overline{M}(0,n)$ and a $+1$ to $\overline{M}(-1,n)$ and $\overline{M}(1,n)$.
They present the example of the overpartition $7+\overline{5}+\overline{2}+1$ which contributes a $-1$ to the count of $\overline{M}(0,15)$ and a $+1$ to $\overline{M}(-1,15)$ and $\overline{M}(1,15)$. We make this adjustment more explicit with a further example.
Table~\ref{t2} lists the overpartitions of $3$ excluding $\overline{2}+1$ and the cranks of their subpartitions consisting of the non-overlined parts. We then evaluate the first residual crank obtained after making the appropriate adjustment explained above. More precisely, the overpartition $\overline{2}+1$ contributes a~$-1$ to the count of $\overline{M}(0,3)$ and a $+1$ to $\overline{M}(-1,3)$ and $\overline{M}(1,3)$.
Then
\begin{alignat*}{4}
&\overline{M}(-3,3)=1,\qquad&&\overline{M}(-2,3)=1,\qquad&&\overline{M}(-1,3)=0,&
\\
&\overline{M}(0,3)=3,&
\\
&\overline{M}(1,3)=1,\qquad&&\overline{M}(2,3)=1,\qquad&&\overline{M}(3,3)=0.&
\end{alignat*}
Making the adjustment described above, we get
\begin{alignat*}{4}
&\overline{M}(-3,3)=1,\qquad&&\overline{M}(-2,3)=1,\qquad&&\overline{M}(-1,3)=1,&
\\
&\overline{M}(0,3)=2,&
\\
&\overline{M}(1,3)=1,\qquad&&\overline{M}(2,3)=1,\qquad&&\overline{M}(3,3)=1,&
\end{alignat*}
which agrees with our calculations and the coefficient of $q^3$ in the expansion of the generating function of the first residual crank.

\begin{table}[t]\centering
\caption{Crank of subpartitions consisting of non-overlined parts of overpartitions of $3$.}\label{t2}
\vspace{2mm}
\begin{tabular}{ c c c }
Overpartitions of $3$ & Crank of the subpartition consisting \\ &of the non-overlined parts \\
 $3$                         & $3$ \\
 $\overline{3}$              & $0$ \\
 $2+1$                       & $0$ \\
 $2+\overline{1}$            & $2$ \\
 $\overline{2}+\overline{1}$ & $0$ \\
 $1+1+1$                     & \hspace{-3mm}$-3$ \\
 $\overline{1}+1+1$          & \hspace{-3mm}$-2$
\end{tabular}
\end{table}

\section{Proof of Theorem \thm{M2}}\label{section3}

\begin{proof}[\unskip\nopunct]
The generating function for the second residual crank of overpartitions is
\begin{align*}
\overline{C2}(z,q)&=\frac{(-q;q)_{\infty}}{\bigl(q;q^2\bigr)_{\infty}}C\bigl(z;q^2\bigr),
\end{align*}
where from equation \eqn{crankexp}, we have
\begin{align*}
&{}C(z,q)=\dfrac{(q;q)_{\infty}}{(zq;q)_{\infty}\bigl(z^{-1}q;q\bigr)_{\infty}}
(1-q+zq)+\sum_{n=1}^\infty \sum_{m=-n}^n  M(m,n) z^m q^n.
\end{align*}
Then
\begin{align*}
\overline{C2}(z,q)&=(-q;q)_{\infty} \biggl(\bigl(1-q^2+zq^2\bigr)+\sum\limits_{\substack{{\pi \in \mathcal{P}}\\ \pi \neq \varnothing}}z^{\operatorname{crank}(\pi)}q^{2 |\pi|}\biggr) \sum_{\pi \in \mathcal{PO}}q^{|\pi|}
\\
&=(-q;q)_{\infty} \bigl(1-q^2+zq^2\bigr)\sum_{\pi \in \mathcal{PO}}q^{|\pi|}+(-q;q)_{\infty} \sum\limits_{\substack{{\pi \in \mathcal{PE}}\\ \pi \neq \varnothing}}z^{\operatorname{crank}(\frac{\pi}{2})}q^{|\pi|} \sum_{\pi \in PO}q^{|\pi|},
\end{align*}
where the partition $\frac{\pi}{2}$ is obtained by dividing each part of the partition $\pi \in \mathcal{PE}$ by $2$. We look for a combinatorial interpretation of $(-q;q)_{\infty}\bigl(1-q^2+zq^2\bigr)$.

Removing the first two parts from equation \eqn{distexp}, we get
\begin{align}
\mylabel{eq:RFP2}
&\sum\limits_{m=3}^{\infty}(1+q)\cdots\bigl(1+q^{m-1}\bigr)q^m=(-q;q)_{\infty}-1-q^2(1+q)-q.
\end{align}
Also, we have the following identity:
\begin{align}
\mylabel{eq:EID2}
&\sum\limits_{m=3}^{\infty}(1+q)\cdots\bigl(1+q^{m-3}\bigr)q^m=q^2\sum\limits_{m=1}^{\infty}(1+q)\cdots\bigl(1+q^{m-1}\bigr)q^{m}=q^2((-q;q)_{\infty}-1).
\end{align}
Then, using equations \eqn{RFP2} and \eqn{EID2}, we get
\begin{align*}
(-q;q)_{\infty}\bigl(1-q^2+zq^2\bigr)\!&{}=1+zq^2(-q;q)_{\infty}+(-q;q)_{\infty}-1-q^2-q^2((-q;q)_{\infty}-1)
\\
&{}=1+q+q^3+zq^2((-q;q)_{\infty}-1+1)
\\
&\quad{}+((-q;q)_{\infty}-1-q^2(1+q)-q)-q^2((-q;q)_{\infty}-1)
\\
&{}=1+q+q^3+zq^2
\\
&\quad{}+\!\sum\limits_{m=3}^{\infty}(1+q)\cdots\bigl(1+q^{m-3}\bigr)q^{m}\bigl(z+\bigl(1+q^{m-2}\bigr)\bigl(1+q^{m-1}\bigr)-1\bigr)
\\
&{}=1+q+zq^2+q^{1+2}+zq^3
\\
&\quad{}+\!\sum\limits_{m=4}^{\infty}(1+q)\cdots\bigl(1+q^{m-3}\bigr)q^{m}\bigl(z+q^{m-1}\!+q^{m-2}\!+q^{(m-1)+(m-2)}\bigr)
\\
&{}=\!\sum\limits_{\pi \in \mathcal{PD}}z^{\kappa(\pi)}q^{|\pi|}.
\end{align*}
We proceed as in the proof of Theorem \thm{M1} except we write an overpartition as an element $\vec{\pi}=(\pi_1,\pi_2,\pi_3) \in \mathcal{PD} \times \mathcal{PE} \times \mathcal{PO}$.
Thus,
\begin{align*}
\overline{C2}(z,q)&{}=\sum_{n=0}^\infty
\sum_{m=-n}^n \overline{M2}(m,n) z^m q^n
\\
&{}=(-q;q)_{\infty}\bigl(1-q^2+zq^2\bigr)\frac{1}{\bigl(q;q^2\bigr)_{\infty}}+(-q;q)_{\infty}\sum_{n=1}^\infty \sum_{m=-n}^n M(m,n) z^m q^{2n}\frac{1}{\bigl(q;q^2\bigr)_{\infty}}
\\
&{}=\sum\limits_{\pi_1 \in \mathcal{PD}}z^{\kappa(\pi_1)}q^{|\pi_1|}\sum\limits_{\pi_3 \in \mathcal{PO}}q^{|\pi_3|}+\sum\limits_{\pi_1 \in \mathcal{PD}}q^{|\pi_1|}\sum\limits_{\substack{\pi_2 \in \mathcal{PE}\\ \pi_2 \neq \varnothing}}z^{\operatorname{crank}(\frac{\pi_2}{2})}q^{|\pi_2|}\sum\limits_{\pi_3 \in \mathcal{PO}}q^{|\pi_3|}
\\
&{}=\sum\limits_{\vec{\pi} = (\pi_1,\pi_2,\pi_3) \in \mathcal{PD} \times \mathcal{PE} \times \mathcal{PO}}z^{\operatorname{crank2}(\vec{\pi})}q^{|\vec{\pi}|},
\end{align*}
where equation \eqn{crank2} is used, and for a partition $\pi$ into distinct parts, $\kappa(\pi)$ is defined into three subcases as in Definition \refdef{crank2}. Theorem \thm{M2} follows.
\end{proof}

\begin{table}[hbt!]\centering
\caption{Second residual crank of overpartitions of $4$.}\label{t3}
\vspace{2mm}
\begin{tabular}{ c c c c c c }
Overpartitions of $4\,(\vec{\pi})$ & $\pi_1$ & $\kappa(\pi_1)$ & $\pi_2$ & $\operatorname{crank}(\frac{\pi_2}{2})$ & $\operatorname{crank2}(\vec{\pi})$ \\
 $4$                          & $\varnothing$ &     & $4$    & $2$ & $2$\\
 $\overline{4}$               & $4$    & $1$ & $\varnothing$ &     & $1$\\
 $3+1$                        & $\varnothing$ &     & $\varnothing$ &     & $0$\\
 $\overline{3}+1$             & $3$    & $1$ & $\varnothing$ &     & $1$\\
 $3+\overline{1}$             & $1$    & $0$ & $\varnothing$ &     & $0$\\
 $\overline{3}+\overline{1}$  & $3+1$  & $0$ & $\varnothing$ &     & $0$\\
 $2+2$                        & $\varnothing$ &     & $2+2$  & $-2$& \hspace{-3mm}$-2$\\
 $\overline{2}+2$             & $2$    & $1$ & $2$    & $-1$& \hspace{-3mm}$-1$\\
 $2+1+1$                      & $\varnothing$ &     & $2$    & $-1$& \hspace{-3mm}$-1$\\
 $\overline{2}+1+1$           & $2$    & $1$ & $\varnothing$ &     & $1$\\
 $2+\overline{1}+1$           & $1$    & $0$ & $2$    & $-1$& \hspace{-3mm}$-1$\\
 $\overline{2}+\overline{1}+1$& $2+1$  & $0$ & $\varnothing$ &     & $0$\\
 $1+1+1+1$                    & $\varnothing$ &     & $\varnothing$ &     & $0$\\
 $\overline{1}+1+1+1$         & $1$    & $0$ & $\varnothing$ &     & $0$\\\\
\end{tabular}
\end{table}

We present an example.
Expanding the generating function of the second residual crank, we find that the coefficient of $q^4$ is
\begin{equation*}
\frac{z^4+3z^3+6z^2+3z+1}{z^2},
\end{equation*}
which agrees with the crank values found in Table \ref{t3}.

\section[M\_2-rank and crank of partitions without repeated odd parts]{$\boldsymbol{M_2}$-rank and crank of partitions without repeated odd parts}\label{section4}

In 2008, Lovejoy and Osburn initiated a series of papers on finding formulas for generating functions of partition statistics' differences using the approach developed by Atkin and Swinnerton-Dyer \cite{At-Sw}. In one of these papers \cite{Lo-Os09}, they proved formulas for the generating functions for $M_2$-rank differences for partitions without repeated odd parts. The $M_2$-rank for partitions without repeated odd parts was first defined by Berkovich and the first author \cite{Be-Ga}. The $M_2$-rank of such partitions is defined to be the number of columns minus the number of rows of its $2$-modular diagram and if $N2(m, n)$ denotes the number of partitions of $n$ without repeated odd parts whose $M_2$-rank is $m$, then the generating function is given by \cite[equation~(1.1)]{Lo-Os09}
\begin{equation*}
N2(z,q) = \sum_{n=0}^\infty
\sum_{m=-n}^n N2(m,n) z^m q^n=\sum\limits_{n=0}^{\infty}q^{n^2}\dfrac{\bigl(-q;q^2\bigr)_{n}}{\bigl(zq^2;q^2\bigr)_{n}\bigl(z^{-1}q^2;q^2\bigr)_{n}}.
\end{equation*}
Then, a natural $M_2$-crank of partitions without repeated odd parts (say M2crank) is given by
\begin{equation*}
M2(z,q) = \dfrac{\bigl(q^2;q^2\bigr)_{\infty}}{\bigl(zq^2;q^2\bigr)_{\infty}\bigl(z^{-1}q^2;q^2\bigr)_{\infty}}\bigl(-q;q^2\bigr)_{\infty}.
\end{equation*}
This crank is not new and appears in the work of the first author and Jennings-Shaffer \cite{Ga-Sh}, where they define the $M_2$-crank as follows. For a partition $\pi$ of $n$ without repeated odd parts
we take the crank of the partition $\frac{\pi_e}{2}$ obtained
by taking the subpartition $\pi_e$, of the even parts of $\pi$,
and halving each part of $\pi_e$. Then, if $M2(m,n)$ is the number
of partitions $\pi$ of $n$ without repeated odd parts such that the
partition $\frac{\pi_e}{2}$ has crank $m$, then the generating
function two variable generating function for $M2(m,n)$ is the
formula above. However they note that this interpretation is not quite correct and it fails for partitions with distinct odd parts whose only even parts are a single two. Here, we deduce the exact combinatorial interpretation of the M2crank.

We define \begin{equation*}
M2(z,q) = \dfrac{\bigl(q^2;q^2\bigr)_{\infty}}{\bigl(zq^2;q^2\bigr)_{\infty}\bigl(z^{-1}q^2;q^2\bigr)_{\infty}}\bigl(-q;q^2\bigr)_{\infty} = (z-1)q^2 + \sum_{\substack{n=0\\ n\neq 2}}^\infty \sum_{m=-n}^n M2(m,n) z^m q^n.
\end{equation*}

\subsection{Connections with tenth order mock theta functions}

In \cite[p.~9]{Ra88}, Ramanujan gave a list of identities involving four tenth order mock theta functions~$\phi(q)$,~$\psi(q)$,~$X(q)$ and $\chi(q)$. In Section~\ref{section1}, we saw that the rank and the first residual crank of overpartitions had connections with the tenth order mock theta functions $\phi(q)$ and $\psi(q)$. In~this section, we consider the remaining two tenth order mock theta functions and establish their connections with the $M_2$-rank of partitions without repeated odd parts. If $N2(a,m,n)$ and $M2(a,m,n)$ denote the number of partitions of $n$ without repeated odd parts whose $M_2$-rank and respectively $M_2$-crank are congruent to $a$ modulo $m$, then one can deduce that
\begin{align*}
\sum_{n\geq 0}(N2(0,5,5n)-N2(2,5,5n))q^n&=X(-q),
\end{align*}
where
\begin{equation*}
X(q)=\sum\limits_{n\geq 0}\frac{(-1)^nq^{n^2}}{(-q;q)_{2n}}
\end{equation*}
is a tenth order mock theta function of Ramanujan. This follows from the $M_2$-rank difference identity of Lovejoy and Osburn \cite[equation~(1.11)]{Lo-Os09} and the identity
\begin{equation*}
\chi(q)=2-qb_2(q)-2q^2k\bigl(q^2,q^5\bigr)
\end{equation*}
for the tenth order mock theta function $X(q)$ due to Choi ({see \cite[p.~183]{Ch00} for details}). Similarly, one can deduce that
\begin{align*}
\sum_{n\geq 0}(N2(1,5,5n+4)-N2(2,5,5n+4))q^n=q^{-1}\chi(-q),
\end{align*}
where
\[
\chi(q)=\sum\limits_{n\geq 0}\frac{(-1)^nq^{(n+1)^2}}{(-q;q)_{2n+1}},
\]
is the other tenth order mock theta function of Ramanujan. Again, this follows from the $M_2$-rank difference identity of Lovejoy and Osburn \cite[equation~(1.10)]{Lo-Os09} and the identity
\begin{equation*}
\chi(q)=2-qb_2(q)-2q^2k\bigl(q^2,q^5\bigr)
\end{equation*}
for the tenth order mock theta function $\chi(q)$ due to Choi (see \cite[p.~183]{Ch00} for details).

The Jacobi products that appear in the formulas for $M_2$-rank differences for partitions without repeated odd parts due to Lovejoy and Osburn \cite[Theorem~1.2]{Lo-Os09} also appear in formulas for $M_2$-crank differences for partitions without repeated odd parts due to the first author and Jennings-Shaffer \cite[Theorem 2.12]{Ga-Sh}. Combining these identities and the above two identities due to Choi, we can deduce relations between the $M_2$-rank and crank of partitions without repeated odd parts and the tenth order mock theta functions $X(q)$ and $\chi(q)$:
\begin{align*}
&{}2\sum_{n\geq 0}(N2(0,5,5n)-N2(1,5,5n))q^n+\sum_{n\geq 0}(M2(0,5,5n)-M2(1,5,5n))q^n=3X(-q),
\\
&{}2\sum_{n\geq 0}(N2(0,5,5n+4)-N2(1,5,5n+4))q^n
\\
&\qquad{}-\sum_{n\geq 0}(M2(0,5,5n+4)-M2(1,5,5n+4))q^n=q^{-1}\chi(-q).\\
\end{align*}

\subsection[Combinatorial interpretation of the M2crank (M\_2-crank for partitions without repeated odd parts)]{Combinatorial interpretation of the M2crank \\($\boldsymbol{M_2}$-crank for partitions without repeated odd parts)}

From the previous sections, we can deduce that
\begin{align*}
M2(z,q)&=\bigl(-q;q^2\bigr)_{\infty} C\bigl(z;q^2\bigr)
\\
&=(-q;q^2)_{\infty} \bigl(1-q^2+zq^2\bigr)+\bigl(-q;q^2\bigr)_{\infty} \sum\limits_{\substack{{\pi \in \mathcal{PE}}\\ \pi \neq \varnothing}}z^{\operatorname{crank}(\frac{\pi}{2})}q^{|\pi|},
\end{align*}
where the partition $\frac{\pi}{2}$ is obtained by dividing each part of the partition $\pi \in \mathcal{PE}$ by $2$. We look for a combinatorial interpretation of $\bigl(-q;q^2\bigr)_{\infty}\bigl(1-q^2+zq^2\bigr)$. We know
\begin{equation*}
\bigl(-q;q^2\bigr)_{\infty}=1+\sum\limits_{m=1}^{\infty}(1+q)\bigl(1+q^3\bigr)\cdots\bigl(1+q^{2m-3}\bigr)q^{2m-1},
\end{equation*}
where $2m-1$ corresponds to the largest part of a partition into distinct odd parts. Therefore,
\begin{align*}
&\sum\limits_{m=1}^{\infty}(1+q)\bigl(1+q^3\bigr)\cdots\bigl(1+q^{2m-3}\bigr)q^{2m-1}=\bigl(-q;q^2\bigr)_{\infty}-1.
\end{align*}
Removing the first part, we get
\begin{align}
\mylabel{eq:RFP3}
&\sum\limits_{m=2}^{\infty}(1+q)\bigl(1+q^3\bigr)\cdots\bigl(1+q^{2m-3}\bigr)q^{2m-1}=\bigl(-q;q^2\bigr)_{\infty}-1-q.
\end{align}
Also, we have the following identity:
\begin{align}
\nonumber
\sum\limits_{m=2}^{\infty}(1+q)\cdots\bigl(1+q^{2m-5}\bigr)q^{2m-1}&=\sum\limits_{m=1}^{\infty}(1+q)\cdots\bigl(1+q^{2m-3}\bigr)q^{2m+1}
\\
\nonumber
&=q^2\sum\limits_{m=1}^{\infty}(1+q)\cdots\bigl(1+q^{2m-3}\bigr)q^{2m-1}
\\
\mylabel{eq:EID3}
&=q^2\bigl(\bigl(-q;q^2\bigr)_{\infty}-1\bigr).
\end{align}
Then, using equations \eqn{RFP3} and \eqn{EID3}, we get
\begin{align*}
(-q;q^2)_{\infty}\bigl(1-q^2+zq^2\bigr)&=1+zq^2\bigl(-q;q^2\bigr)_{\infty}+\bigl(-q;q^2\bigr)_{\infty}-1-q^2-q^2\bigl(\bigl(-q;q^2\bigr)_{\infty}-1\bigr)
\\
&=1+q-q^2+zq^2\bigl(\bigl(-q;q^2\bigr)_{\infty}-1+1\bigr)
\\
&\quad+\bigl(\bigl(-q;q^2\bigr)_{\infty}-1-q)-q^2\bigl(\bigl(-q;q^2\bigr)_{\infty}-1\bigr)
\\
&=1+q+(z-1)q^2
\\
&\hspace{5mm}+\sum\limits_{m=2}^{\infty}(1+q)\bigl(1+q^3\bigr)\cdots\bigl(1+q^{2m-5}\bigr)q^{2m-1}\bigl(z+q^{2m-3}\bigr)
\\
&=1+q+(z-1)q^2+(z+q)q^3
\\
&\quad+\sum\limits_{m=3}^{\infty}(1+q)\bigl(1+q^3\bigr)\cdots\bigl(1+q^{2m-5}\bigr)q^{2m-1}\bigl(z+q^{2m-3}\bigr)
\\
&=1+q+(z-1)q^2+(z+q)q^3
\\
&\quad+\sum\limits_{m=1}^{\infty}(1+q)\bigl(1+q^3\bigr)\cdots\bigl(1+q^{2m-1}\bigr)q^{2m+3}\bigl(z+q^{2m+1}\bigr)
\\
&=(z-1)q^2+\sum\limits_{\pi \in \mathcal{POD}}z^{\theta(\pi)}q^{|\pi|},
\end{align*}
where $\mathcal{POD}$ denotes the set of partitions into distinct odd parts, and $\theta(\pi)$ for $\pi \in \mathcal{POD}$  is defined into the following two subcases.

Case I: $\ell(\pi) \geq 5$,
\begin{equation*}
\theta(\pi)=\begin{cases}
0 & \mbox{if $\ell(\pi)-2$ is a part of $\pi$}, \\
1 & \mbox{otherwise.}
\end{cases}
\end{equation*}

Case II: $\ell(\pi) = 3$ or $1$,
\begin{equation*}
\theta(\pi)=\begin{cases}
0 & \mbox{if $1$ is a part}, \\
1 & \mbox{otherwise.}
\end{cases}
\end{equation*}
We consider a partition without repeated odd parts as an element
$\vec{\pi}=(\pi_1,\pi_2) \in \mathcal{POD} \times \mathcal{PE}$.
Then
\begin{align*}
M2(z,q)&{}=(z-1)q^2+\sum_{\substack{n=0\\ n\neq 2}}^\infty
\sum_{m=-n}^n  M2(m,n) z^m q^n
\\
&{}=(z-1)q^2+\sum\limits_{\pi_1 \in \mathcal{POD}}z^{\theta(\pi_1)}q^{|\pi_1|}+\sum\limits_{\pi_1 \in \mathcal{POD}}q^{|\pi_1|}\sum\limits_{\substack{\pi_2 \in \mathcal{PE}\\ \pi_2 \neq \varnothing}}z^{\operatorname{crank}(\frac{\pi_2}{2})}q^{|\pi_2|}
\\
&{}=(z-1)q^2+\sum\limits_{\vec{\pi} = (\pi_1,\pi_2) \in \mathcal{POD} \times \mathcal{PE}}z^{\operatorname{M2crank}(\vec{\pi})}q^{|\vec{\pi}|},
\end{align*}
where
\begin{equation*}
\operatorname{M2crank}(\vec{\pi})=\begin{cases}
\operatorname{crank}\bigl(\frac{\pi_2}{2}\bigr) & \mbox{if $\pi_2 \neq \varnothing$},
\\
\theta(\pi_1) & \mbox{if $\pi_2 = \varnothing$.}
\end{cases}
\end{equation*}
The following theorem follows easily.
\begin{Theorem}
The number of partitions $\vec{\pi}$ of $n$ without repeated odd parts with $\operatorname{M2crank}(\vec{\pi})=m$ is $M2(m,n)$ for all $n \neq 2$.
\end{Theorem}

\section{Conclusion}

In this paper, we have seen that various rank and their companion crank statistics have explicit relations to mock theta functions of order $5$ and $10$. A problem worth considering is to explore more results of this type and relate other mock theta functions to rank and crank differences. Furthermore, explicit $p$-dissections of rank-type statistics discussed in this paper for higher primes $p \geq 7$ are scarce in the literature and might be another problem worth pursuing.

\subsection*{Acknowledgements}
We thank the referees for their comments and suggestions.

\pdfbookmark[1]{References}{ref}
\LastPageEnding

\end{document}